\documentclass[12pt,amsfonts]{amsart}
\usepackage{graphicx,amssymb,eucal}
\usepackage{latexsym,amsmath,epsfig,epic,eepic}
\usepackage[all]{xy}
\usepackage[english, french]{babel}

\newtheorem{theorem}{{\bf Th\'eor\`eme}}[section]
\newtheorem{main}{{\bf  Th\'eor\`eme}}
\newtheorem{proposition}[theorem]{{\bf Proposition}}

\newcommand{\RR}{\mathbb{R}^{n}}
\newcommand{\C}{\mathbb{C}}
\newcommand{\R}{\mathbb{R}}
\newcommand{\Z}{\mathbb{Z}}
\newcommand{\SP}{\mathrm{Sp}(2n,\mathbb{R})}

\newcommand{\ang}{\mathrm{angle}}
\newcommand{\ham}{\mathrm{Ham}}
\newcommand{\calabi}{\mathfrak{Cal}}
\newcommand{\symp}{\mathrm{Symp}}

\newcommand{\SPE}{\mathrm{Sp(E,\omega)}}


\title{Quasi-morphismes et invariant de Calabi}
\author{Pierre Py}
\date{Juin 2005}

\begin{document}
\maketitle

\selectlanguage{english}
\begin{abstract} 
In this paper, we give two elementary constructions of homogeneous quasi-morphisms defined on the group of Hamiltonian diffeomorphisms of certain closed connected symplectic manifolds (or on its universal cover). The first quasi-morphism, denoted by $\calabi_{S}$, is defined on the group of Hamiltonian diffeomorphisms of a closed oriented surface $S$ of genus greater than $1$. This construction is motivated by a question of M. Entov and L. Polterovich \cite{entpol}. If $U\subset S$ is a disk or an annulus, the restriction of $\calabi_{S}$ to the subgroup of diffeomorphisms which are the time one map of a Hamiltonian isotopy in $U$ equals Calabi's homomorphism. The second quasi-morphism is defined on the universal cover of the group of Hamiltonian diffeomorphisms of a symplectic manifold for which the cohomology class of the symplectic form is a multiple of the first Chern class.     
\end{abstract}

\selectlanguage{french}
\begin{abstract}
Dans ce texte, nous donnons deux constructions \'el\'ementaires de quasi-morphismes homog\`enes d\'efinis sur le groupe des diff\'eomorphismes hamiltoniens d'une vari\'et\'e symplectique connexe ferm\'ee (ou sur son rev\^etement universel). Le premier quasi-morphisme, not\'e $\calabi_{S}$, est d\'efini sur le groupe des diff\'eomorphismes hamiltoniens d'une surface ferm\'ee orient\'ee de genre sup\'erieur ou \'egal \`a $2$. Cette construction est motiv\'ee par une question de M. Entov et L. Polterovich \cite{entpol}. Si $U\subset S$ est un disque ou un anneau, la restriction de $\calabi_{S}$ au groupe des diff\'eomorphismes qui sont le temps $1$ d'une isotopie hamiltonienne dans $U$ est \'egale au morphisme de Calabi. Le second quasi-morphisme est d\'efini sur le rev\^etement universel du groupe des diff\'eomorphismes hamiltoniens d'une vari\'et\'e symplectique pour laquelle la classe de cohomologie de la forme symplectique est un multiple de la premi\`ere classe de Chern.   
\end{abstract}

\section{Introduction}
Un \textit{quasi-morphisme} sur un groupe $\Gamma$ est une fonction $\phi : \Gamma \to \R$ telle que les quantit\'es 
$$\phi(xy)-\phi(x)-\phi(y)$$
soient born\'ees lorsque $x,y$ d\'ecrivent $\Gamma$. On appellera parfois \textit{d\'efaut} de $\phi$ la quantit\'e 
$$\delta = \mathrm{sup}_{x,y\in \Gamma}\vert\phi(xy)-\phi(x)-\phi(y)\vert .$$ Un quasi-morphisme est \textit{homog\`ene} s'il satisfait en outre $\phi(x^{n})=n\phi(x)$ pour $x\in \Gamma$ et $n\in \Z$. Nous dirons que deux quasi-morphismes sont \textit{\'equivalents} si leur diff\'erence est born\'ee. Il n'est pas difficile de v\'erifier que si $\phi$ est un quasi-morphisme quelconque, la formule 
$$\phi_{h}(x)=\mathrm{lim}_{p\to \infty}\frac{1}{p}\phi(x^{p})$$ d\'efinit l'unique quasi-morphisme homog\`ene \`a distance born\'ee de $\phi$. Il satisfait $$\vert \phi -\phi_{h}\vert \le \delta.$$ On pourra consulter (par exemple) \cite{bav,brooks} pour une introduction \`a ce sujet. Nous noterons $QM_{h}(\Gamma,\R)$ l'espace des quasi-morphismes homog\`enes sur le groupe $\Gamma$.

Si $(V,\omega)$ est une vari\'et\'e symplectique connexe ferm\'ee, le groupe $\ham(V,\omega)$ de ses diff\'eomorphismes hamiltoniens est simple d'apr\`es un th\'eor\`eme de A. Banyaga \cite{ban1}. Il n'admet donc pas de morphisme non-trivial vers $\R$. Si $(V,\omega)$ est une vari\'et\'e symplectique connexe ouverte, sur laquelle $\omega$ est exacte, E. Calabi a introduit dans \cite{cal} un morphisme 
$$\calabi_{V} : \ham(V,\omega)\to \R. $$ Le noyau de ce morphisme est simple d'apr\`es un autre th\'eor\`eme de Banyaga \cite{ban1}. Si $\lambda$ est une primitive de $\omega$ sur $V$ et $(f_{t})$ une isotopie hamiltonienne dans $V$, engendr\'ee par le champ de vecteurs $Z_{t}$, on a : 
$$\calabi_{V}(f_{1})=\int_V \! \int_{0}^{1} \lambda(Z_{t})dt\, \omega .$$

\noindent Supposons maintenant que $(V,\omega)$ est ferm\'ee. A chaque ouvert connexe $U\subset V$, on associe le sous-groupe $\Gamma_{U}$ de $\ham(V,\omega)$ form\'e des diff\'eomorphismes qui sont le temps $1$ d'une isotopie hamiltonienne dans $U$. Si $\omega$ est exacte sur $U$, on dispose alors du morphisme de Calabi : $\calabi_{U} : \Gamma_{U} \to \R$. On notera $\mathcal{D}$ la famille des ouverts connexes $U$ de $V$ tels que $\omega$ est exacte sur $U$ et tels qu'il existe $f\in \ham(V,\omega)$ avec $f(U)\cap \overline{U}=\emptyset$. Dans \cite{entpol}, M. Entov et L. Polterovich posent la question suivante : 

\begin{center}
\textit{Peut-on construire un quasi-morphisme homog\`ene $\phi : \ham(V,\omega)\to \R$ dont les restrictions aux sous-groupes $(\Gamma_{U})_{U\in \mathcal{D}}$ soient \'egales aux morphismes de Calabi $(\calabi_{U})_{U\in \mathcal{D}}$ ? }
\end{center}
Plus g\'en\'eralement, peut-on construire un tel quasi-morphisme sur le rev\^etement universel $\widetilde{\ham}(V,\omega)$ du groupe des diff\'eomorphismes hamiltoniens ? Dans \cite{entpol} (voir aussi \cite{biran}), ils r\'epondent positivement \`a cette question pour une certaine classe de vari\'et\'es symplectiques, qui inclut notamment les espaces projectifs complexes, en particulier la sph\`ere $S^{2}$. Leur m\'ethode utilise des outils sophistiqu\'es de topologie symplectique. 

Dans l'esprit de constructions pr\'ec\'edentes de J.-M. Gambaudo et \'E. Ghys, nous construisons un quasi-morphisme ayant une propri\'et\'e comparable \`a celle \'enonc\'ee dans la question ci-dessus, sur le groupe des diff\'eomorphismes hamiltoniens d'une surface de genre sup\'erieur ou \'egal \`a $2$. Dans \cite{gg1}, on pourra trouver de nombreuses constructions de quasi-morphismes sur les groupes $\ham(S,\omega)$, pour toute surface ferm\'ee orient\'ee (les quasi-morphismes de Gambaudo et Ghys sont en fait d\'efinis sur le groupe $\symp_{0}(S,\omega)$).  Cependant les restrictions de ces quasi-morphismes aux sous-groupes $(\Gamma_{U})$ ne sont pas des homomorphismes. Les constructions de Gambaudo et Ghys peuvent \^etre vues comme de possibles g\'en\'eralisations du nombre de translation $\tau$ sur le groupe $\widetilde{\mathrm{Hom\acute{e}o}}_{+}(S^{1})$ des hom\'eomorphismes de la droite r\'eelle qui commutent aux translations enti\`eres.  Elles sont dans l'esprit de constructions pr\'ec\'edentes par V.I. Arnold \cite{arnold2} ,  D. Ruelle \cite{ruelle}, S. Schwartzman \cite{schw}.

\begin{main} Soit $S$ une surface ferm\'ee orient\'ee de genre sup\'erieur ou \'egal \`a $2$, munie d'une forme symplectique $\omega$. Il existe un quasi-morphisme homog\`ene 
$$\calabi_{S} : \ham(S,\omega)\to \R,$$ dont la restriction aux sous-groupes $\Gamma_{U}$ est \'egale au morphisme de Calabi, d\`es que $U$ est diff\'eomorphe \`a un disque ou \`a un anneau. Le quasi-morphisme $\calabi_{S}$ est invariant par conjugaison par tout diff\'eomorphisme symplectique. 
\end{main}

Nous pouvons faire deux remarques concernant ce th\'eor\`eme. D'une part, en utilisant des constructions de Gambaudo et Ghys, on peut le renforcer en l'\'enonc\'e :

\textit{l'espace (affine) des quasi-morphismes homog\`enes sur $\ham(S,\omega)$ ayant la propri\'et\'e du th\'eor\`eme $1$ est de dimension infinie.} 

\noindent Cela r\'esultera simplement du fait que l'espace des quasi-morphismes homog\`enes sur $\ham(S,\omega)$ dont les restrictions aux groupes $\Gamma_{U}$ (o\`u $U$ est diff\'eomorphe \`a un disque ou \`a un anneau) sont nulles, est un espace de dimension infinie. D'autre part, la nature de ce quasi-morphisme est certainement tr\`es diff\'erente de celle du quasi-morphisme construit par Entov et Polterovich dans \cite{entpol}. En effet, parmi les deux conditions qui d\'efinissent la famille d'ouverts $\mathcal{D}$, la premi\`ere est vide en dimension $2$, en revanche la seconde n'appara\^it pas du tout dans notre travail. Il existe d'ailleurs des disques (ou des anneaux) ne la satisfaisant pas. 

L'\'enonc\'e du th\'eor\`eme suivant  est inspir\'e du th\'eor\`eme $5.2$ de \cite{entpol}, o\`u un calcul similaire est fait sur le groupe $\ham(S^{2},\omega)$. Notre m\'ethode est cependant diff\'erente.  Consid\'erons donc  une fonction de Morse $F : S \to \R$, dont les valeurs critiques sont toutes distinctes. Notons $x_{1},\ldots, x_{l},$ ses points critiques, $\lambda_{j}=F(x_{j})$ ses valeurs critiques, avec $\lambda_{1}< \cdots < \lambda_{l}$. On a classiquement : 
$$\sum_{j=1}^{l}(-1)^{\mathrm{ind}x_{j}}=2-2g,$$ o\`u $g$ est le genre de $S$. Consid\'erons l'espace 
$$\mathcal{F}=\{H : S\to \R, \{H,F\}=0\},$$ des fonctions sur la surface $S$ qui commutent avec $F$ au sens de Poisson, c'est-\`a-dire telles que $\omega (X_{H},X_{F})=0$,
o\`u $X_{G}$ d\'esigne le gradient symplectique d'une fonction $G$. L'ensemble 
$$\Gamma=\{\varphi_{H}^{1}, H \in \mathcal{F}\}$$
est un sous-groupe  ab\'elien de $\ham(S,\omega)$ (o\`u $\varphi_{H}^{t}$ d\'esigne le flot de $X_{H}$). La restriction de $\calabi_{S}$ \`a $\Gamma$ est donc un homomorphisme, que nous calculons dans le th\'eor\`eme suivant. La donn\'ee de la fonction $F$ permet de construire un graphe fini $\mathcal{G}$ appel\'e graphe de Reeb \cite{reeb}, de la mani\`ere suivante. Parmi les composantes connexes des niveaux $F^{-1}(cste)$ on trouve : 

\begin{enumerate}
\item les points critiques de $F$ d'indice $0$ ou $2$,
\item des courbes simples plong\'ees,
\item des courbes immerg\'ees ayant un unique point double (correspondant \`a un point critique d'indice $1$ de $F$).
\end{enumerate}
A chaque composante de type $1$ ou $3$ on associe un sommet de $\mathcal{G}$. Notons $K$ la r\'eunion des composantes de type $1$ ou $3$. L'ouvert $S\setminus K$ est une r\'eunion finie de cylindres diff\'eomorphes \`a $S^{1}\times \R$. A chaque cylindre $C$ on associe une ar\^ete dont les extr\'emit\'es sont les sommets associ\'es aux composantes de niveaux de $F$ qui contiennent $\partial C$. Nous avons une application naturelle $p_{\mathcal{G}} : S \to \mathcal{G}$, et, si $H\in \mathcal{F}$, on peut \'ecrire $H=H_{\mathcal{G}}\circ p_{\mathcal{G}}$, o\`u $H_{\mathcal{G}}$ est d\'efinie sur $\mathcal{G}$.

Nous pouvons \textit{\'elaguer} le graphe $\mathcal{G}$ pour obtenir un graphe $\mathcal{G}'$, de la mani\`ere suivante. Le graphe $\mathcal{G}$ poss\`ede des sommets de degr\'e $1$ ou $3$. Si $v$ est un sommet de degr\'e $1$ de $\mathcal{G}$, nous retirons $v$ ainsi que l'ar\^ete \`a laquelle il \'etait reli\'e, pour obtenir un nouveau graphe. Ce faisant, nous pouvons cr\'eer un sommet de degr\'e $2$ (ou un autre sommet de degr\'e $1$ \`a partir de la seconde it\'eration de ce proc\'ed\'e). R\'ep\'etons ce proc\'ed\'e jusqu'\`a obtenir un graphe $\mathcal{G}'$ qui ne poss\`ede plus que des points de degr\'e $2$ ou $3$. Les sommets de degr\'e $3$ de $\mathcal{G}'$ sont en nombre $2g-2$. En effet la quantit\'e 
$$\sum_{v}2-degr\acute{e} (v),$$
o\`u la somme porte sur tous les sommets de $\mathcal{G}$, est \'egale \`a la caract\'eristique d'Euler de la surface, et reste constante au cours de l'\'elagage. On note $\mathcal{V}$ l'ensemble de ces $2g-2$ sommets. Nous supposons dans le th\'eor\`eme suivant que l'aire totale de la forme $\omega$ est \'egale \`a $2g-2$.

\begin{main} Si $H$ est dans $\mathcal{F}$, nous avons :
$$\calabi_{S}(\varphi_{H}^{1})=\int_{S}H\omega -\sum_{v\in \mathcal{V}}H_{\mathcal{G}}(v).$$
\end{main}

Dans la derni\`ere  partie, nous proposons une construction \'el\'ementaire  d'un quasi-morphisme homog\`ene d\'efini sur le rev\^etement universel du groupe des diff\'eo\-mor\-phismes hamiltoniens d'une vari\'et\'e symplectique connexe ferm\'ee $(V,\omega)$ pour laquelle $[\omega] = r c_{1}(V)$.  Ici, $[\omega]$ d\'esigne la classe de cohomologie de la forme symplectique et $c_{1}(V)$ la premi\`ere classe de Chern de $V$. Pour les besoins de la cause, nous \'ecrivons notre hypoth\`ese sous la forme : 
$$s[\omega]=2c_{1}(V),$$ o\`u $s$ est un r\'eel non-nul. Nous calculons la restriction de ce quasi-morphisme sur les isotopies hamiltoniennes support\'ees dans une boule ; son expression fait alors intervenir un quasi-morphisme $\tau_{B,\omega}$ introduit par J. Barge et \'E. Ghys, dont la construction est rappel\'ee dans le texte. Si $(f_{t})$ est une isotopie hamiltonienne dans $V$, on note $\{f_{t}\}$ l'\'el\'ement du rev\^etement universel du groupe $\ham(V,\omega)$ qu'elle d\'efinit.

\begin{main}Si $(V,\omega)$ v\'erifie l'hypoth\`ese ci-dessus, il existe un quasi-morphisme homog\`ene 
$$\mathfrak{S} : \widetilde{\ham}(V,\omega)\to \R,$$ tel que, pour toute isotopie hamiltonienne $\{f_{t}\}$ support\'ee dans une boule $B$, on ait : $$\mathfrak{S}(\{f_{t}\})=\tau_{B,\omega}(f_{1})+s\calabi_{B}(f_{1}).$$
\end{main}
 
Le texte est organis\'e comme suit. Dans la seconde partie, nous prouvons les th\'eor\`emes $1$ et $2$, apr\`es avoir rappel\'e des r\'esultats de Banyaga. Dans la troisi\`eme partie, nous rappelons des constructions de Barge et Ghys, puis nous prouvons le th\'eor\`eme $3$. Exception faite des r\'esultats mentionn\'es dans le paragraphe $2.1$, qui sont utilis\'es dans la troisi\`eme partie, les parties $2$ et $3$ sont ind\'ependantes.

Une annonce des th\'eor\`emes $1$ et $2$ est contenue dans \cite{py}.

\section{Surfaces de genre sup\'erieur}
\subsection{Extension du groupe des diff\'eomorphismes hamiltoniens}

Les r\'esultats de ce paragraphe sont dus \`a Banyaga \cite{ban2} ; nous les rappelons succinctement.

Consid\'erons une vari\'et\'e (ferm\'ee, connexe) $M$ munie d'une forme de contact $\alpha$ dont le champ de Reeb $X$ est induit par une action libre du cercle $\R/\Z$. La vari\'et\'e $V$, quotient de $M$ par l'action du cercle, porte une forme symplectique $\omega$ telle que $\pi^{*}\omega=d\alpha$, o\`u $\pi : M\to V$ est la projection canonique.

Dans cette situation nous avons une extension centrale par $\R/\Z$ du groupe des diff\'eomorphismes hamiltoniens de $V$. D\'ecrivons d'abord cette extension au niveau des alg\`ebres de Lie de champs de vecteurs. Un \'el\'ement $Y$ de l'alg\`ebre  $\mathcal{L}_{\alpha}(M)$ des champs de vecteurs sur $M$ qui pr\'eservent $\alpha$ est invariant par l'action du cercle $\R/\Z$, il d\'efinit un champ de vecteurs hamiltonien $\pi_{*}Y$ sur $V$. L'alg\`ebre $\mathcal{L}_{\alpha}(M)$ est donc une extension centrale par $\R$ de l'alg\`ebre $\mathrm{ham}(V,\omega)$ des champs de vecteurs hamiltoniens sur $V$. Si $Z\in \mathrm{ham}(V,\omega)$ (avec $\iota_{Z}\omega=dH_{Z}$, $\int_{V}H_{Z}\omega^{n}=0$, $\mathrm{dim}V=2n$), le champ de vecteurs $\theta(Z)=\widehat{Z}-(H_{Z}\circ \pi)X$  (o\`u $\widehat{Z}$ est le relev\'e horizontal de $Z$ : $\alpha(\widehat{Z})=0$) pr\'eserve $\alpha$. L'application $Z\mapsto \theta(Z)$ est un morphisme d'alg\`ebres de Lie qui scinde l'extension.

Notons $G_{\alpha}(M)_{0}$ le groupe des diff\'eomorphismes de $M$ qui pr\'eservent $\alpha$, isotopes \`a l'identit\'e via une isotopie qui pr\'eserve $\alpha$. On a alors une extension centrale : 
\centerline{
\xymatrix{ {0} \ar[r]Ê& {\R/\Z} \ar[r] & {G_{\alpha}(M)_{0}}Ê\ar[r] & {\ham (V,\omega)} \ar[r] & {0}. \\
}
}

\noindent Si $(f_{t})$ est une isotopie hamiltonienne on note $\Theta(f_{t})$ l'isotopie de $M$ obtenue en ``int\'egrant'' $\theta$. Sa classe d'homotopie ne d\'epend que de celle de $(f_{t})$. Ainsi, si le groupe $\ham(V,\omega)$ est simplement connexe, l'extension pr\'ec\'edente est scind\'ee. Gr\^ace au th\'eor\`eme de Banyaga qui assure que le groupe $\ham(V,\omega)$ est simple, la section qui scinde l'extension est unique. Si $V$ est une surface de genre sup\'erieur, le groupe $\ham(V,\omega)$ est simplement connexe \cite{earle,mcduff}, et l'extension ci-dessus est canoniquement scind\'ee.

\subsection{Construction du quasi-morphisme $\calabi_{S}$}

On suppose donc que $S$ est une surface ferm\'ee orient\'ee, de genre $g$ sup\'erieur ou \'egal \`a $2$, munie d'une forme symplectique d'aire totale $2g-2$ (d'apr\`es un th\'eor\`eme de Moser deux telles formes sont l'image l'une de l'autre par un diff\'eomorphisme isotope \`a l'identit\'e, le choix de $\omega$ est donc sans importance).

Nous noterons $M$ la vari\'et\'e des droites orient\'ees tangentes \`a $S$, $\widetilde{S}$ le rev\^etement universel de $S$, et $\widetilde{M}$ la vari\'et\'e des droites orient\'ees tangentes \`a $\widetilde{S}$. Le choix d'une m\'etrique \`a courbure constante sur $S$, de forme d'aire \'egale \`a $\omega$, fournit une structure de $S^{1}$-fibr\'e principal sur $M$ et $\widetilde{M}$. Nous noterons $X$ le champ de vecteurs sur $M$ tangent aux fibres de l'application $\pi : M \to S$ engendr\'e par cette action du cercle. On note \'egalement $S^{1}_{\infty}$ le cercle \`a l'infini de $\widetilde{S}$ d\'etermin\'e par cette m\'etrique, et $p_{\infty} : \widetilde{M}\to S^{1}_{\infty}$ la projection. La vari\'et\'e $\widetilde{M}$ est diff\'eomorphe \`a $\widetilde{S}\times S^{1}_{\infty}$ et le feuilletage $(\widetilde{S}\times \{\ast \})_{\ast \in S^{1}_{\infty}}$ descend en un feuilletage sur $M$ appel\'e feuilletage horocyclique. Nous noterons $\alpha$ une $1$-forme sur $M$ telle que $\alpha(X)=1$ et $\pi^{*}\omega=d\alpha$. La forme $\alpha$ est une forme de contact de champ de Reeb \'egal \`a $X$. On notera encore $X$ le relev\'e \`a $\widetilde{M}$ de ce champ. Enfin, si $\gamma : [0,1]\to S^{1}_{\infty}$ est un chemin continu, nous noterons $n(\gamma)$ l'entier d\'efini comme suit. Si un param\'etrage de $S^{1}_{\infty}$ par $\R/\Z$ est donn\'e, notons $\widetilde{\gamma}$ un relev\'e de $\gamma$ \`a $\R$. On pose : 
  $$n(\gamma)=[ \widetilde{\gamma}(1)-\widetilde{\gamma}(0)] ,$$
  cet entier ne d\'epend pas du choix du param\'etrage. Si $\gamma$ et $\beta$ sont deux chemins dans $S^{1}_{\infty}$ avec $\gamma(1)=\beta(0)$ nous avons :
  $$(\ast) \; \; \; \vert n(\gamma \ast \beta)-n(\gamma)-n(\beta)\vert \le 2.$$

\textbf{Observations}. Avant de construire le quasi-morphisme annonc\'e dans le th\'eor\`eme $1$, nous commen\c{c}ons par quelques remarques. A chaque ouvert connexe $U\subset S$, distinct de $S$, nous allons associer un \'el\'ement canonique de l'espace 
$$\frac{QM_{h}(\pi_{1}(U),\R)}{\mathrm{Hom}(\pi_{1}(U),\R)}.$$ Puisqu'un quasi-morphisme homog\`ene est invariant par conjugaison, nous oublierons parfois de choisir un point base pour le groupe $\pi_{1}(U)$. Soit $\psi : U\times S^{1}\to \pi^{-1}(U)$ une trivialisation du fibr\'e $\pi : M\to S$ au-dessus de $U$. Si $z_{0}$ est un point sur $S^{1}$ et si $\gamma$ est un lacet dans $U$, bas\'e en $x_{0}$, on pose : $$\phi_{z_{0}}([\gamma])=n(p_{\infty}(\widetilde{\psi(\gamma(t),z_{0})})),$$
  o\`u $\widetilde{\psi(\gamma(t),z_{0})}$ est un relev\'e \`a $\widetilde{M}$ du chemin $\psi (\gamma(t),z_{0})$. Vu la propri\'et\'e $(\ast)$, l'application $\phi_{z_{0}}$ est un quasi-morphisme sur le groupe $\pi_{1}(U,x_{0})$. Son homog\'en\'eis\'e $\phi$ ne d\'epend pas du point $z_{0}$. Si deux m\'etriques \`a courbure constante sur $S$ sont donn\'ees, on peut identifier par un hom\'eomorphisme les bords \`a l'infini de $\widetilde{S}$ associ\'es, l'indice $n$ ne d\'epend donc pas de la m\'etrique. A l'addition d'un homomorphisme pr\`es, il ne d\'epend pas de la trivialisation $\psi$, la classe 
  $$[\phi] \in \frac{QM_{h}(\pi_{1}(U),\R)}{\mathrm{Hom}(\pi_{1}(U),\R)}$$ est donc canonique. Une autre mani\`ere de la d\'ecrire serait la suivante. On fixe un point $\widetilde{x_{0}}\in M$ au-dessus de $x_{0}$. Si $[\gamma]\in \pi_{1}(U,x_{0})$, on note $\widetilde{\gamma}$ son relev\'e issu de $\widetilde{x_{0}}$ tangent au feuilletage horocyclique. On peut \'ecrire $\widetilde{\gamma}(t)=\psi(\gamma(t),z(t))$. Notons $f([\gamma])$ la variation de l'argument de $z(t)$ compt\'ee en tours. Alors $f$ est un quasi-morphisme dont l'homog\'en\'eis\'e est \'egal \`a $-\phi$. On peut en quelque sorte penser \`a  un feuilletage transverse aux fibres du fibr\'e $M\to S$ comme \`a une ``quasi-trivialisation"  du fibr\'e. Alors que la comparaison de deux trivialisations du fibr\'e au-dessus de l'ouvert $U$ fournit un homomorphisme $\pi_{1}(U)\to \Z$, la comparaison du feuilletage avec une trivialisation fournit un quasi-morphisme homog\`ene sur le groupe $\pi_{1}(U)$. Si le groupe $\pi_{1}(U)$ est ab\'elien, il n'admet pas de quasi-morphisme homog\`ene non-trivial (c'est-\`a-dire autre que les homomorphismes). En revanche si le groupe $\pi_{1}(U)$ est libre non-ab\'elien, la classe $[\phi]$ va appara\^itre comme une obstruction \`a ce que le quasi-morphisme $\calabi_{S} : \ham(S,\omega) \to \R$ que nous allons d\'efinir se restreigne en le morphisme de Calabi sur le groupe $\Gamma_{U}$.  

Donnons enfin une derni\`ere interpr\'etation de la classe $[\phi]$. Rappelons qu'\`a toute repr\'esentation $\rho$ d'un groupe discret $\Gamma$ dans le groupe $\mathrm{Hom\acute{e}o}_{+}(S^{1})$ des hom\'eomor\-phismes du cercle qui pr\'eservent l'orientation, on peut associer une classe de cohomologie born\'ee $eu_{b}(\rho) \in H_{b}^{2}(\Gamma,\Z)$, appel\'ee classe d'Euler born\'ee, voir \cite{ghys}. La repr\'esentation $\pi_{1}(S)\to \mathrm{PSL}(2,\R)$ associ\'ee \`a une m\'etrique \`a courbure constante sur $S$ fournit donc une classe $eu_{b}(S)\in H_{b}^{2}(\pi_{1}(S), \Z)$, qui ne d\'epend pas de la m\'etrique. Pour tout groupe discret $\Gamma$, le noyau de l'application $H_{b}^{2}(\Gamma,\R)\to H^{2}(\Gamma,\R)$ est isomorphe \`a l'espace $$\frac{QM_{h}(\Gamma,\R)}{\mathrm{Hom(\Gamma,\R)}}.$$ Si $U\subset S$ est un ouvert connexe distinct de $S$, on note $i_{U} : \pi_{1}(U)\to \pi_{1}(S)$ le morphisme naturel. Puisque le groupe $H^{2}(\pi_{1}(U),\R)$ est trivial, la classe $i_{U}^{*}eu_{b}(S)$ (que nous consid\'erons comme une classe r\'eelle) est dans le noyau 
$$\mathrm{Ker}(H_{b}^{2}(\pi_{1}(U),\R)\to H^{2}(\pi_{1}(U),\R)).$$
C'est la classe $[\phi]$ pr\'ec\'edemment d\'ecrite.

\textit{Nous pouvons maintenant construire le quasi-morphisme $\calabi_{S}$}. Soit $(f_{t})$ une isotopie hamiltonienne dans $S$. Notons $\Theta(f_{t}) : M\to M$ l'isotopie qui rel\`eve $(f_{t})$ pr\'ec\'edemment construite, et $(F_{t})$ l'isotopie de $\widetilde{M}$ qui rel\`eve $\Theta(f_{t})$. Si $v$ et $w$ sont deux points de $\widetilde{M}$ tels que $\widetilde{\pi}(v)=\widetilde{\pi}(w)$ (o\`u $\widetilde{\pi}$ est la projection $\widetilde{M} \to \widetilde{S}$), on \'ecrit $w=\phi_{X}^{u_{0}}(v)$ avec $u_{0}\in [0,1]$, o\`u $\phi_{X}^{u}$ est le flot de $X$. Notons $G(u,t)=p_{\infty}(F_{t}(\phi_{X}^{uu_{0}}v))$. Le lacet lu sur le ``bord" de $G$ a un indice $n$ \'egal \`a $0$. Puisque les chemins $G(-,0)$ et $G(-,1)$ ont un indice born\'e par $1$, on a $\vert n(p_{\infty}(F_{t}(v)))-n(p_{\infty}(F_{t}(w))\vert \le 2$. On d\'efinit alors, pour $\widetilde{x}\in \widetilde{S}$, $\widetilde{\ang}(\widetilde{x},f_{1})=-\mathrm{inf}_{\widetilde{\pi}(v)=\widetilde{x}}n(p_{\infty}(F_{t}(v)))$. La fonction $\widetilde{\ang}(-,f_{1})$ est invariante sous l'action du groupe fondamental de $S$ et d\'efinit une fonction mesurable born\'ee $\ang(-,f_{1})$ sur $S$. Elle v\'erifie :
$$\vert \ang(x,f_{1}g_{1})-\ang(x,g_{1})-\ang(g_{1}(x),f_{1})\vert \le 8. $$
Ainsi l'application qui, au diff\'eomorphisme hamiltonien $f_{1}$ associe l'int\'egrale $$\int_{S}\ang(-,f_{1})\, \omega$$ est un quasi-morphisme. Nous pouvons l'homog\'en\'eiser pour d\'efinir : 
$$\calabi_{S}(f_{1})=\mathrm{lim}_{p\to \infty}\frac{1}{p}\int_{S}\ang(-,f_{1}^{p})\, \omega. $$ D'apr\`es le th\'eor\`eme ergodique sous-additif \cite{kingman,oseledec}, la suite de fonctions 
$$\frac{1}{p}\ang(-,f^{p})$$ converge $\omega$-presque partout quand $p$ tend vers l'infini, vers une fonction mesurable $\widehat{\ang}(-,f)$. Il n'est pas difficile de v\'erifier que cette fonction est born\'ee et satisfait : 
$$\calabi_{S}(f)=\int_{S}\widehat{ \ang}(-,f)\, \omega .$$  

Discutons maintenant des diff\'erents choix effectu\'es pour notre construction. Nous avons d\'ej\`a indiqu\'e pourquoi l'indice $n$ d'une courbe ne d\'epend pas du choix de la m\'etrique. Par ailleurs :\begin{itemize}
\item Si l'on change de m\'etrique, l'action du cercle sur $M$ (i.e. le champ $X$)  change, mais la classe d'Euler du fibr\'e n'\'etant pas modifi\'ee, on peut trouver un diff\'eomor\-phisme de $M$, induisant l'identit\'e sur $S$, qui entrelace les deux actions. On en d\'eduit ais\'ement l'invariance du quasi-morphisme.
\item Lorsque la m\'etrique est fix\'ee, le choix de la forme $\alpha$ est sans importance. Une autre primitive de $\pi^{*}\omega$ valant $1$ sur $X$ serait de la forme $\alpha + \pi^{*}\beta$, o\`u $\beta$ est une $1$-forme ferm\'ee sur la surface. En utilisant la nullit\'e du flux d'une isotopie hamiltonienne, on voit que le quasi-morphisme final est inchang\'e.  Notons que, lorsque $M$ est identifi\'e au fibr\'e unitaire tangent \`a $S$, la forme $\alpha$  peut-\^etre prise nulle dans la direction du flot g\'eod\'esique. 
\item Une fois acquise l'ind\'ependance de $\calabi_{S}$ vis-\`a-vis de la m\'etrique, l'invariance par conjugaison dans $\symp(S,\omega)$ est claire. Il suffit de consid\'erer une m\'etrique (\`a courbure constante, de forme d'aire $\omega$) et de la transporter par le diff\'eomor\-phisme symplectique consid\'er\'e. 
\end{itemize}

\smallskip

Supposons que $U\subset S$ soit un ouvert connexe distinct de $S$ et $(f_{t})$ une isotopie hamiltonienne dans $U$. Nous allons calculer $\calabi_{S}(f_{1})$.

Choisissons une trivialisation $\psi : U\times S^{1}\to \pi^{-1}(U)$ du fibr\'e $M$ au-dessus de $U$, telle que $X=\frac{\partial}{\partial s}$ (o\`u $s$ est la coordonn\'ee angulaire sur le cercle). On a alors \mbox{$\alpha =ds +\pi^{*}\lambda$}, o\`u $\lambda$ est une primitive de $\omega$ sur $U$. Notons \'egalement $\phi$ le quasi-morphisme homog\`ene sur $\pi_{1}(U)$ associ\'e \`a cette trivialisation. Enfin, on note $Z_{t}$ le champ de vecteurs qui engendre l'isotopie, $H_{t}$ un hamiltonien pour $Z_{t}$ avec $\mathrm{supp}(H_{t})\subset U$ et $\widetilde{H}_{t}$ la fonction d'int\'egrale nulle sur $S$ qui diff\`ere de $H_{t}$ par une constante. 

Fixons un point $x_{0}$ dans $U$ et un compact $K$ tel que $\mathrm{supp}(f_{t})\subset K$. Pour tout point $x$ de $K$, on choisit un chemin $\alpha_{x_{0}x}$ de classe $C^{1}$ par morceaux de $x_{0}$ \`a $x$, contenu dans $U$, de d\'eriv\'ee born\'ee ind\'ependamment de $x$.  On note $\gamma_{x,f}$ le lacet $\alpha_{x_{0}x}\ast (f_{t}(x))\ast \overline{\alpha_{x_{0}f(x)}}$ et $\langle [\phi],f\rangle (x)$ la limite de la suite $(\frac{1}{p}\phi([\gamma_{x,f^{p}}]))_{p\ge 0}$ (pour relier $f^{p}$ \`a l'identit\'e on utilise bien s\^ur l'isotopie $(f_{t})$ concat\'en\'ee $p$ fois). 
\begin{proposition} Pour presque tout $x$ de $U$ nous avons : 
$$\widehat{\ang}(x,f)=-\langle [\phi],f \rangle (x)+\mathfrak{M}(y\mapsto \int_{0}^{1}(\lambda(Z_{t})+\widetilde{H}_{t})(f_{t}(y))dt)(x).$$
\end{proposition}
Dans cette proposition, et dans la suite, $\mathfrak{M}(\varphi)$ d\'esignera la limite des moyennes de Birkhoff d'une fonction (int\'egrable) $\varphi$ relativement \`a la transformation $f=f_{1}$. 

\textit{Preuve} : Notons $(h_{t})$ l'isotopie obtenue en concat\'enant $p$ fois $(f_{t})$. Si $x$ est dans $U$ et $z$ dans $S^{1}$ nous avons :

$$\Theta(h_{t})(\psi(x,z))=\psi(h_{t}(x),exp(-2i\pi \int_{0}^{t}(\lambda(Z_{t'})+\widetilde{H}_{t'})(h_{t'}(x))dt')\cdot z).$$ 

\noindent Notons $v(t)$, $v _{1}(t)$, $v_{2}(t)$ les courbes $$\Theta(h_{t})(\psi(x,z)),$$ $$\psi(h_{t}(x),z),$$ $$\psi(f^{p}(x),exp(-2i\pi \int_{0}^{t}(\lambda(Z_{t'})+\widetilde{H}_{t'})(h_{t'}(x))dt')),$$ respectivement. On choisit deux relev\'es $\widetilde{v}_{1}$ et $\widetilde{v}_{2}$ \`a $\widetilde{M}$ tels que $\widetilde{v}_{1}(1)=\widetilde{v}_{2}(0)$, et l'on note $\widetilde{v}=\widetilde{v}_{1}\ast \widetilde{v}_{2}$. Dans la suite d'\'egalit\'es ci-dessous, le symbole $\simeq$ voudra dire que les deux membres de l'\'egalit\'e diff\`erent d'une quantit\'e born\'ee, la valeur exacte de la borne important peu (mais pouvant cependant ais\'ement \^etre d\'etermin\'ee). 
$$\ang(x,f^{p})\simeq-n(p_{\infty}(\widetilde{v}(t)))\simeq -n(p_{\infty}(\widetilde{v}_{1}(t)))-n(p_{\infty}(\widetilde{v}_{2}(t))),$$
$$-n(p_{\infty}(\widetilde{v}_{2}(t)))\simeq \int_{0}^{1}(\lambda(Z_{t'})+\widetilde{H}_{t'})(h_{t'}(x))dt',$$
$$-n(p_{\infty}(\widetilde{v}_{1}(t)))\simeq -\phi([\gamma_{x,f^{p}}]).$$ Nous obtenons au total l'existence d'une constante $C$ telle que : 
$$\vert \ang(x,f^{p})-\int_{0}^{1}(\lambda(Z_{t'})+\widetilde{H}_{t'})(h_{t'}(x))dt'+\phi([\gamma_{x,f^{p}}])\vert \le C.$$ Le r\'esultat suit. 
\hfill $\Box$

Une fois la proposition pr\'ec\'edente acquise, nous pouvons conclure la preuve du th\'eor\`eme $1$.  Hors de $U$, la fonction $\widehat{\ang}(x,f)$ est \'egale \`a $\int_{0}^{1}\widetilde{H}_{t}(x)dt$. Nous obtenons donc : 
$$\calabi_{S}(f)=-\int_{U}\langle [\phi],f \rangle \omega +\calabi_{U} (f)$$ (la fonction $\langle [\phi],f \rangle$ est bien s\^ur nulle hors de l'ouvert $U$). Si $U$ est simplement connexe, le terme $\langle [\phi],f\rangle$ est identiquement nul. Si $U$ est un anneau, le quasi-morphisme homog\`ene $\phi$ est un morphisme qui se repr\'esente par une $1$-forme ferm\'ee $\beta$. Dans ce cas, on v\'erifie sans peine que  $$\langle [\phi],f \rangle(x)=\mathfrak{M}(y \mapsto \int_{0}^{1}\beta (Z_{t})(f_{t}(y))dt)(x).$$ On en d\'eduit : 
 $$\int_{U}\langle [\phi],f\rangle \omega = \int_{U}\int_{0}^{1}\beta(Z_{t})\, dt \, \omega.$$ Mais cette derni\`ere int\'egrale est nulle pour une isotopie hamiltonienne. Dans le cas d'un flot hamiltonien, cela traduit simplement le fait que le cycle asymptotique de Schwartzman pour la mesure $\omega$ est nul. On a donc bien $\calabi_{S}(f_{1})=\calabi_{U}(f_{1})$ dans ce cas. Si le groupe $\pi_{1}(U)$ est libre non-ab\'elien, l'espace 
 $$\frac{QM_{h}(\pi_{1}(U),\R)}{\mathrm{Hom}(\pi_{1}(U),\R)}$$ n'est pas trivial \cite{bav,brooks}, et l'int\'egrale $\int_{U}\langle [\phi,f]\rangle \omega$ peut ne pas s'annuler. Nous avons achev\'e la preuve du th\'eor\`eme $1$.
 
 D\'ecrivons maintenant les constructions de Gambaudo et Ghys qui permettent de prouver que l'espace des quasi-morphismes homog\`enes sur $\ham(S,\omega)$ nuls en restriction aux sous-groupes $(\Gamma_{U})$ o\`u $U$  est un disque ou un anneau, est de dimension infinie.  
 Consid\'erons une $1$-forme $\eta$, non n\'ecessairement ferm\'ee, sur la surface $S$, et notons $\widetilde{\eta}$ son rel\`evement \`a $\widetilde{S}$. On suppose toujours qu'une m\'etrique \`a courbure constante de forme d'aire $\omega$ est fix\'ee.  Nous allons construire un quasi-morphisme homog\`ene 
$$\Phi_{\eta} : \symp_{0}(S, \omega) \to \R.$$  Un \'el\'ement $f$ de $\symp_{0}(S,\omega)$ admet un relev\'e canonique $\widetilde{f}$ \`a $\widetilde{S}$ : puisque le groupe $\symp_{0}(S,\omega)$ est simplement connexe nous choisissons une isotopie $(f_{t})$ reliant l'identit\'e \`a $f$, et nous consid\'erons son relev\'e $(\widetilde{f}_{t})$ \`a $\widetilde{S}$. On pose alors $\widetilde{f}=\widetilde{f}_{1}$. Si $x$ est dans $\widetilde{S}$ nous noterons $\delta (x,f)$ l'unique g\'eod\'esique (pour la m\'etrique fix\'ee) qui relie $x$ \`a $\widetilde{f}(x)$. La fonction 
$$x\mapsto \int_{\delta(x,f)}\widetilde{\eta},$$ est $\pi_{1}(S)$-invariante et d\'efinit donc une fonction $u(\eta, f)$ sur la surface $S$.
En utilisant le fait que la $2$-forme $d\widetilde{\eta}$ v\'erifie une in\'egalit\'e de la forme $\parallel d\widetilde{\eta} \parallel \le C \parallel \widetilde{\omega}\parallel$ sur $\widetilde{S}$, et le fait que les triangles g\'eod\'esiques de $\widetilde{S}$ sont d'aire born\'ee, nous obtenons ais\'ement :
$$\vert u(\eta,fg)-u(\eta,g)-u(\eta,f)\circ g \vert \le \pi C.$$
Ainsi $f\mapsto \int_{S}u(\eta,f)\omega$ d\'efinit un quasi-morphisme sur $\symp_{0}(S, \omega)$, son homog\'en\'eis\'e sera not\'e $\Phi_{\eta}$ (contrairement \`a $\calabi_{S}$, il d\'epend de la m\'etrique).  

Dans \cite{gg1}, il est montr\'e  que la famille $(\Phi_{\eta})_{\eta}$, restreinte au groupe $\ham(S,\omega)$, engendre un espace de dimension infinie dans l'espace $QM_{h}(\ham(S,\omega),\R)$. Les auteurs utilisent pour cela des ``twists" support\'es au voisinage d'une g\'eod\'esique ferm\'ee simple. Notons que, pour que ceux-ci soient hamiltoniens, cette g\'eod\'esique doit \^etre homologue \`a $0$.

Supposons que $(f_{t})$ est une isotopie \textit{symplectique} support\'ee dans le compact $K$ contenu dans l'ouvert simplement connexe $U$ de $S$. Notons $\widetilde{K}\subset \widetilde{U}$ des relev\'es \`a $\widetilde{S}$. Le compact $\widetilde{K}$ est stable par $\widetilde{f}$ et de diam\`etre fini. On a donc $\vert u(\eta,f^{p}) \vert \le \vert \eta \vert \cdot \mathrm{diam}(\widetilde{K})$ pour tout $p$. Ceci assure que $\Phi_{\eta}(f_{1})=0$. Soit maintenant $(f_{t})$ une isotopie \textit{hamiltonienne} (engendr\'ee par le champ $Z_{t}$)  support\'ee dans l'anneau $A=]0,1[\times \R/\Z \hookrightarrow S$. On suppose bien s\^ur ce plongement injectif au niveau du groupe fondamental, sans quoi on serait ramen\'e au cas pr\'ec\'edent. On v\'erifie alors que $$\widehat{u}(\eta,f)=\mathrm{lim}_{p\to \infty}\frac{1}{p}u(\eta,f^{p})$$ vaut $l\cdot \mathfrak{M}(\int_{0}^{1}\beta(Z_{t})\circ f_{t} dt)$, o\`u la classe de la $1$-forme ferm\'ee $\beta$ sur $A$ engendre $H^{1}(A,\Z)$ et $l$ est l'int\'egrale de la $1$-forme $\eta$ sur la g\'eod\'esique ferm\'ee librement homotope dans $S$ au g\'en\'erateur de $\pi_{1}(A)$. On en d\'eduit $\Phi_{\eta}(f_{1})=0$.

\subsection{Calcul sur des hamiltoniens autonomes}

Nous prouvons ici le th\'eor\`eme $2$. Notons $U$ l'ouvert $S-\{x_{l}\}$ et fixons une trivialisation du fibr\'e $\pi : M\to S$ au-dessus de $U$. Celle-ci fournit une primitive $\lambda$ de $\omega$ sur $U$ et un quasi-morphisme homog\`ene $\phi$ sur le groupe $\pi_{1}(U)$, comme pr\'ec\'edemment. Si $a$ est une ar\^ete du graphe $\mathcal{G}$ nous noterons $a^{+}$ et $a^{-}$ les sommets \`a ses extr\'emit\'es (avec la convention $F_{\mathcal{G}}(a^{-})<F_{\mathcal{G}}(a^{+})$). Pour chaque ar\^ete $a$, nous fixons un param\'etrage de $p_{\mathcal{G}}^{-1}(a)$ par $\R/\Z \times ]F_{\mathcal{G}}(a^{-}),F_{\mathcal{G}}(a^{+})[$, avec des coordonn\'ees $(\theta,t)$ telles que $F(\theta,t)=t$ et $\omega= d\theta \wedge dt$. Si $H$ est une fonction dans $\mathcal{F}$, le champ hamiltonien $Z_{H}$ s'\'ecrit sur $p_{\mathcal{G}}^{-1}(a)$ : 
$$\vartheta(t)\frac{\partial}{\partial \theta},$$ o\`u la fonction $\vartheta$ satisfait :
$$\int_{F_{\mathcal{G}}(a^{-})}^{F_{\mathcal{G}}(a^{+})}\vartheta(t)dt=H_{\mathcal{G}}(a^{+})-H_{\mathcal{G}}(a^{-}).$$

Pour un domaine \`a bord lisse $D\subset U$ nous noterons $\langle [\phi], \partial D \rangle$ la somme des valeurs de $\phi$ sur les classes de conjugaison d\'etermin\'ees par chacune des composantes du bord de $D$. Cette valeur ne d\'epend que de la classe $[\phi]$.

\begin{proposition} Si $D$ est \`a bord g\'eod\'esique, pour une m\'etrique \`a courbure constante quelconque sur $S$, on a $\langle [\phi], \partial D\rangle =-\chi(D)$.
\end{proposition}

\textit{Preuve} : on peut supposer que $\omega$ est la forme d'aire associ\'ee \`a la m\'etrique qui rend le bord de $D$ g\'eod\'esique (car le quasi-morphisme $\phi$ est ind\'ependant de $\omega$). On a vu que la $1$-forme $\alpha$ peut alors \^etre choisie nulle dans la direction du flot g\'eod\'esique. Au-dessus de $U$, dans une trivialisation dans laquelle $X=\partial / \partial s$, on a $\alpha = ds+\pi^{*}(\lambda)$. Soit $\gamma $ une orbite p\'eriodique du flot g\'eod\'esique telle que $\pi(\gamma)$ est une composante du bord de $D$. Puisque $\gamma$ est ferm\'ee on a : $\phi([\pi(\gamma)])=-\int_{\gamma}ds=\int_{\pi(\gamma)}\lambda$. En sommant sur les diff\'erentes composantes de bord, on obtient : $\langle [\phi], \partial D\rangle =\int_{D}\omega =-\chi(D)$. \hfill $\Box$

Si $H$ est dans $\mathcal{F}$, bien que l'isotopie $(\varphi_{H}^{t})$ ne soit pas n\'ecessairement \`a support dans $U$, on peut r\'ep\'eter le raisonnement qui a servi  \`a \'etablir la proposition $2.1$. Si $x$ est dans $U$ et $p_{\mathcal{G}}(x)$ n'est pas un sommet de $\mathcal{G}$, on notera $[x]$ la classe d'homotopie libre du cercle $p_{\mathcal{G}}^{-1}(p_{\mathcal{G}}(x))$, orient\'e par $X_{F}$. Notant $\widetilde{H}$ la fonction d'int\'egrale nulle sur $S$ qui diff\`ere de $H$ par une constante, on a, presque partout sur $U$ : 
$$\widehat{\ang}(x,\varphi_{H}^{1})=\mathfrak{M}(y\mapsto \lambda(X_{H})(y)+\widetilde{H}(y))(x)-\vartheta (x)\phi([x]).$$ La fonction $\mathfrak{M}(y\mapsto \lambda(X_{H})(y)+\widetilde{H}(y))$ a pour int\'egrale $\int_{S}H\omega-(2g-2)H(x_{l})$. Si $a$ est une ar\^ete de $\mathcal{G}$, nous noterons $[a]$ la valeur commune des classes $[x]$ pour $x\in p_{\mathcal{G}}^{-1}(a)$. La somme des int\'egrales de $\vartheta (x)\phi([x])$ sur les diff\'erentes ar\^etes vaut : 
$$\sum_{a}\phi([a])(H_{\mathcal{G}}(a^{+})-H_{\mathcal{G}}(a^{-})).$$ On peut \'ecrire la somme sous la forme $$\sum_{v}C(v)H_{\mathcal{G}}(v),$$ o\`u la somme porte cette fois sur les sommets de $\mathcal{G}$. Il faut calculer les constantes $C(v)$. Au pr\'ealable, notons le fait suivant :

\textbf{Observation}. Si $x$ et $y$ sont deux lacets dans $\pi^{-1}(U)$, tangents au feuilletage horocyclique de $S$, ayant m\^eme point base, alors $\phi([\pi(x\ast y)])=\phi([\pi(x)])+\phi([\pi(y)])$.

\noindent On peut alors commencer le calcul des constantes $C(v)$ :
\begin{itemize}
\item Si $v$ correspond \`a un extremum local autre que le maximum global, la constante $C(v)$ est nulle. En effet, elle est \'egale \`a la valeur de $\phi$ sur un petit lacet qui entoure l'extremum. Comme celui-ci est dans $U$ le lacet est nul dans $\pi_{1}(U)$.
\item Si $v=p_{\mathcal{G}}(x_{l})$, la constante $C(v)$ est \'egale \`a la valeur de $\phi$ sur la classe d'homotopie dans $U$ d'un petit lacet qui entoure $x_{l}$ (avec l'orientation oppos\'ee \`a celle du bord de $\{F \le \lambda_{l}-\epsilon\}$). Si ce lacet $\gamma$ est assez petit, il admet un relev\'e tangent au feuilletage horocyclique qui est ferm\'e. On constate alors que $\phi([\gamma])=-(2g-2)$ (la classe d'Euler du fibr\'e en cercles au-dessus de $S$).
\end{itemize}
Il reste alors \`a calculer les constantes $C(v)$ pour les sommets $v$ correspondant \`a des points critiques d'indice $1$. Il n'est pas difficile de v\'erifier que $$C(p_{\mathcal{G}}(x_{j}))=\langle [\phi], \partial \{\lambda_{j}-\epsilon \le F \le \lambda_{j}+ \epsilon \}\rangle,$$ (pour $\epsilon$ assez petit). Le domaine $\{\lambda_{j}-\epsilon \le F \le \lambda_{j}+ \epsilon \}$ est constitu\'e d'un nombre fini de cylindres sur le bord desquels $\phi$ est nul et d'un pantalon $P$.  Il faut \'evaluer le terme $\langle [\phi],\partial P \rangle $.

\begin{itemize}
\item Si l'une des composantes du bord de $P$ est homotope \`a $0$ dans $U$, $\phi$ est nul \'evalu\'e contre celle-ci. Les deux autres composantes sont alors librement homotopes dans $U$, et les deux valeurs de $\phi$ correspondantes sont oppos\'ees. On peut donc supposer les trois composantes de $\partial P$ essentielles dans $U$ sans quoi $C(v)=0$.   
\item Si en outre ces trois composantes sont essentielles dans $S$, on peut trouver une m\'etrique \`a courbure constante (de forme d'aire \'egale \`a $\omega$) qui rend le bord de $P$ g\'eod\'esique (ou seulement deux composantes sur trois de $\partial P$ si deux d'entre elles sont librement homotopes). Une g\'en\'eralisation imm\'ediate de la proposition $2.2$ permet alors de montrer que la constante $C(v)$ vaut $1$.
\item Il reste \`a traiter le cas o\`u l'une des composantes, disons $\alpha_{1}$, de $\partial P$ est contractile dans $S$. Dans ce cas elle borde un disque plong\'e qui contient le point $x_{l}$ (puisque l'on a suppos\'e cette m\^eme courbe essentielle dans $U$).  Les deux autres composantes $\alpha_{2}$ et $\alpha_{3}$ de $\partial P$ sont alors essentielles dans $S$. Modifions le pantalon $P$ en un pantalon $P'$ de composantes de bord $(\alpha_{i}')_{1\le i\le 3}$ telles que $\alpha_{i}'$ est homotope \`a $\alpha_{i}$, et $\alpha_{1}'$ et $\alpha_{2}'$ ont m\^eme point base. Le lacet $\alpha_{1}'$ peut \^etre choisi contenu dans un disque arbitrairement petit au voisinage de $x_{l}$. On peut alors trouver une m\'etrique \`a courbure constante qui rend $\alpha_{2}'$ g\'eod\'esique. Notant $\beta_{2}$ l'orbite p\'eriodique du flot g\'eod\'esique telle que $\pi(\beta_{2})=\alpha_{2}'$, on peut trouver un relev\'e, tangent au feuilletage horocyclique et ferm\'e, $\beta_{1}$ de $\alpha_{1}'$, issu du m\^eme point que l'orbite $\beta_{2}$. Dans ce cas l'observation ci-dessus assure que $\phi([\alpha_{1}'\ast \alpha_{2}'])=\phi([\alpha_{1}'])+\phi([\alpha_{2}'])$. Puisque la derni\`ere composante de bord de $P'$ d\'efinit la classe de conjugaison de $[\alpha_{1}'\ast \alpha_{2}']^{-1}$, on a $\langle [\phi], \partial P' \rangle =0$. 
\end{itemize}
Il n'est pas difficile de v\'erifier que l'ensemble des sommets pour lesquels $C(v)=1$ correspond \`a l'ensemble $\mathcal{V}$ d\'efini dans l'introduction. Une autre mani\`ere de d\'ecrire cet ensemble serait la suivante. Les \'el\'ements de $\mathcal{V}$ sont les sommets correspondant \`a des points critiques d'indice $1$ pour lesquels les trois composantes du bord du pantalon $P$ d\'ecrit pr\'ec\'edemment sont essentielles dans $S$. En r\'esum\'e, la constante $C(p_{\mathcal{G}}(x_{l}))$ vaut $-(2g-2)$, les autres constantes sont \'egales \`a $1$ pour les \'el\'ements de $\mathcal{V}$ et $0$ sinon. 
Finalement, la somme initiale est \'egale \`a : $$\sum_{v\in \mathcal{V}}H_{\mathcal{G}}(v)-(2g-2)H(x_{l}).$$ Nous obtenons ainsi : $$\calabi_{S}(\varphi_{H}^{1})=\int_{S}\widehat{\ang}(-,\varphi_{H}^{1})\omega=\int_{S}H\omega-\sum_{v\in \mathcal{V}}H_{\mathcal{G}}(v).$$

\section{Quasi-morphisme sur certaines vari\'et\'es symplectiques de premi\`ere classe de Chern non-nulle}
\subsection{Exposant de Lyapunov symplectique}

Notons $B$ une boule de $\R^{2n}$, munie d'une forme symplectique $\nu$ de volume fini. Nous rappelons ici une construction due \`a Barge et Ghys \cite{bargeghys}, d'un quasi-morphisme homog\`ene $\tau_{B,\nu}$ sur le groupe $\Gamma_{B,\nu}=\mathrm{Diff}^{c}(B,\nu)$ des diff\'eomorphismes symplectiques de $B$ \`a support compact. Cette construction appara\^it d\'ej\`a implicitement dans \cite{ruelle}. Nous renvoyons \`a \cite{bargeghys} pour plus de d\'etails. 

On commence par construire un quasi-morphisme $\Phi$ sur le rev\^etement universel $\widetilde{\SPE}$ du groupe symplectique d'un espace vectoriel symplectique $(E,\omega)$. Cette construction est bien connue \cite{bargeghys, dupont, guichwigner} ; nous la rappelons cependant, dans un souci de compl\'etude. Supposons que $J$ soit une structure presque-complexe sur $E$ compatible avec $\omega$ : $J$ est un endomorphisme de $E$, de carr\'e $-1$, qui pr\'eserve $\omega$, avec $\omega(u,Ju)>0$ pour tout vecteur non-nul $u$ de $E$. Muni de la forme $(u,v)_{J}=\omega(u,Jv)-i\omega(u,v)$, $E$ devient un espace vectoriel hermitien. Nous noterons $\Lambda(E)$ la grassmannienne lagrangienne de $E$. Si $L_{0}$ et $L_{1}$ sont deux lagrangiens, il existe un endomorphisme unitaire $u$ de $E$  tel que $u(L_{0})=L_{1}$. Le nombre complexe 
$$det_{\C}^{2}(u)\in S^{1},$$ ne d\'epend pas du choix de $u$. On le note $det^{2}_{L_{0}}L_{1}$. Il v\'erifie la relation (de cocycle) : $$det^{2}_{L_{0}}L_{2}=det^{2}_{L_{0}}L_{1}\cdot det^{2}_{L_{1}}L_{2}.$$
En particulier, si $(L_{t})$ est une courbe dans $\Lambda(E)$ la variation de l'argument du nombre complexe $det^{2}_{W}L_{t}$ (compt\'ee en tours) ne d\'epend pas du choix de $W$. On note $\Delta (det^{2} L_{t})$ ce nombre.

\begin{proposition} Si $(L_{t})$ est une courbe dans $\Lambda(E)$ qui reste toujours transverse \`a un lagrangien donn\'e $W$, on a $\vert \Delta(det^{2} L_{t})\vert \le n$.
\end{proposition}
Il est classique que l'application 
$$\begin{array}{rcl}
\Lambda(E) & \to & S^{1} \\
L & \mapsto & det^{2}_{L_{0}}L \\
\end{array}$$
induit un isomorphisme entre les groupes fondamentaux. Par ailleurs, il est \'egalement bien connu que, pour tout lagrangien $W$, l'intersection avec l'hypersurface $$\{L, \mathrm{dim}L\cap W >0\},$$ engendre le groupe $H^{1}(\Lambda(E),\Z)$. Il n'est donc pas surprenant que le fait de rester transverse \`a un lagrangien donn\'e, emp\^eche une courbe de $\Lambda(E)$ de``trop tourner" (voir \cite{arnold}). Prouvons maintenant la proposition. 

\textit{Preuve} :  L'application qui a un endomorphisme $f$ de $W$, sym\'etrique pour le produit scalaire $\langle u,v \rangle_{J}=\omega(u,Jv)$, associe le graphe $L_{f}=\{f(x)+Jx\}_{x\in W}$ est un diff\'eomorphisme sur l'ouvert des lagrangiens transverses \`a $W$. On a :
$$det^{2}_{W}L_{f}=\prod_{k=1}^{n}\frac{(\lambda_{k}+i)^{2}}{1+\lambda_{k}^{2}},$$ o\`u les $\lambda_{k}$ sont les valeurs propres de $f$. Si $(f_{t})$ est un chemin d'endomorphismes sym\'etriques de $W$, de valeurs propres $\lambda_{1}(t)\le \cdots \le \lambda_{n}(t)$, la variation de l'argument du nombre complexe 
$$\prod_{k=1}^{n}\frac{(\lambda_{k}(t)+i)^{2}}{1+\lambda_{k}(t)^{2}}$$  est inf\'erieure \`a $n$. \hfill $\Box$

Un \'el\'ement du rev\^etement universel $\widetilde{\SPE}$ du groupe symplectique $\SPE$ est la donn\'ee d'un chemin $\gamma : [0,1]\to \SPE$, tel que $\gamma(0)=\mathrm{Id}$, d\'efini \`a une homotopie fixant les extr\'emit\'es pr\`es. On le note $[\gamma]\in \widetilde{\SPE}$. On note $\widetilde{\mathrm{Id}}$ l'\'el\'ement d\'efini par le chemin constant \'egal \`a $\mathrm{Id}$. Si $L_{0}\in \Lambda(E)$, on note $\varphi_{L_{0}}([\gamma])=\Delta(det^{2}(\gamma_{t}\cdot L_{0}))$. Si $L_{0}$ et $L_{1}$ sont deux lagrangiens on a 
$$\vert \varphi_{L_{0}}([\gamma])-\varphi_{L_{1}}([\gamma])\vert \le 2n$$ (il suffit de consid\'erer l'application $(s,t)\in[0,1]^{2}\mapsto det^{2}_{\ast}(\gamma_{t}\cdot L_{s})$, o\`u $L_{s}$ est un chemin de $L_{0}$ \`a $L_{1}$ qui reste transverse \`a un lagrangien donn\'e, et d'appliquer la proposition pr\'ec\'edente ). De plus nous avons l'\'egalit\'e 
$$\varphi_{L_{0}}([\gamma]\cdot [\eta])-\varphi_{L_{0}}([\gamma])-\varphi_{L_{0}}([\eta])=\varphi_{\gamma_{1}\cdot L_{0}}([\eta])-\varphi_{L_{0}}([\eta]),$$ cette derni\`ere quantit\'e est born\'ee par $2n$. L'application $\varphi_{L_{0}}$ est donc un quasi-morphisme sur le groupe $\widetilde{\SPE}$, dont l'homog\'en\'eis\'e $\Phi$ ne d\'epend pas de $L_{0}$. Pour tout lagrangien $L$ on a : 
$$\Phi([\gamma])=\mathrm{lim}_{p\to \infty}\frac{1}{p}\Delta(det^{2}(\gamma^{p}_{t}\cdot L)).$$

\begin{proposition}
Le quasi-morphisme $\Phi$ est continu. 
\end{proposition}
\textit{Preuve} : on montre par r\'ecurrence sur $k$ l'in\'egalit\'e : 
$$\vert \varphi_{L_{0}}(x^{kp})-k\varphi_{L_{0}}(x^{p})\vert \le 2nk.$$ En divisant par $kp$ et en faisant tendre $k$ vers l'infini, nous obtenons : 
$$\vert \Phi(x)-\frac{1}{p}\varphi_{L_{0}}(x^{p})\vert \le \frac{2n}{p}.$$ La continuit\'e de $\varphi_{L_{0}}$ implique alors celle de $\Phi$.  \hfill $\Box$

Nous avons vu que le quasi-morphisme $\Phi : \widetilde{\SPE}\to \R$ ne d\'epend d'aucun choix de point base dans la lagrangienne. Il est \'egalement ind\'ependant du choix de la structure presque-complexe. En effet, si $T$ est le g\'en\'erateur du groupe infini cyclique 
$$\pi_{1}(\SPE,\mathrm{Id})\hookrightarrow \widetilde{\SPE},$$ on v\'erifie ais\'ement que $\Phi(T)=2$. Ainsi les deux quasi-morphismes homog\`enes construits \`a partir de deux structures presque-complexes $J$ et $J'$ distinctes, prennent la m\^eme valeur sur l'\'el\'ement $T$. D'apr\`es un argument de Barge et Ghys \cite{bargeghys}, cela entra\^ine qu'ils sont \'egaux. Indiquons finalement que l'on peut donner d'autres descriptions du quasi-morphisme $\Phi$, qui permettent de le calculer effectivement \cite{bargeghys}.

Passons alors \`a la construction du quasi-morphisme homog\`ene $\tau_{B,\nu} : \Gamma_{B,\nu} \to \R$. Nous choisissons une trivialisation symplectique du fibr\'e tangent \`a $B$. Si $f\in \Gamma_{B,\nu}$, la diff\'erentielle de $f$ ``lue" dans cette trivialisation est une application 

$$\begin{array}{rcl}
B & \to & \SP \\
x & \mapsto & df(x). \\
\end{array}$$ 
Un changement de trivialisation est donn\'e par une application $\theta : B \to \SP$. La diff\'erentielle de $f$ est alors chang\'ee en $x \mapsto \theta(f(x))^{-1}\cdot df(x)\cdot \theta(x)$. Notons $\widetilde{\theta} : B \to \widetilde{\SP}$ un relev\'e quelconque de $\theta$, et $\widetilde{df} : B \to \widetilde{\SP}$ l'unique relev\'e de $df$ qui vaut $\widetilde{\mathrm{Id}}$ hors d'un compact. Si $f$ et $g$ sont dans le groupe $\Gamma_{B,\nu}$, nous avons : 
$$\widetilde{d(fg)}(x)=\widetilde{df}(g(x))\cdot \widetilde{dg}(x),$$ et donc 
$$\vert \Phi(\widetilde{d(fg)}(x))-\Phi(\widetilde{dg}(x))-\Phi(\widetilde{df}(g(x)))\vert \le 2n.$$
L'application $f\mapsto \int_{B}\Phi(\widetilde{df})\nu^{n}$ est donc un quasi-morphisme sur le groupe $\Gamma_{B,\nu}$. Nous allons v\'erifier que son homog\'en\'eis\'e ne d\'epend pas de la trivialisation symplectique choisie. Si l'on change de trivialisation, l'application $\widetilde{df}$ est chang\'ee en $\widetilde{\theta}^{-1}\circ f \cdot \widetilde{df} \cdot \widetilde{\theta}$. Nous avons donc, puisque $\Phi$ est homog\`ene : 
$$\vert \Phi(\widetilde{df^{p}}(x)) -\Phi(\widetilde{\theta}^{-1}(f^{p}(x))\cdot \widetilde{df^{p}}(x)\cdot \widetilde{\theta}(x))\vert \le 2n + \vert \Phi(\widetilde{\theta}^{-1}(f^{p}(x))\widetilde{\theta}(x))\vert.  $$ Si le support de $f$ est contenu dans le compact $K$ de $B$, nous avons pour $x$ dans $K$, $\vert \Phi(\widetilde{\theta}^{-1}(f^{p}(x))\widetilde{\theta}(x))\vert \le 2n +2\sup_{K}\vert \Phi\circ \widetilde{\theta}\vert $ ; si $x$ n'est pas dans $K$, la quantit\'e $\Phi(\widetilde{\theta}^{-1}(f^{p}(x))\widetilde{\theta})(x)$ est nulle.  Ainsi : 
$$\vert \Phi(\widetilde{df^{p}}(x)) -\Phi(\widetilde{\theta}^{-1}(f^{p}(x))\widetilde{df^{p}}(x)\widetilde{\theta}(x))\vert \le 4n +2\mathrm{sup}_{K}\vert \Phi\circ \widetilde{\theta}\vert,$$ pour tout $x$ de $B$ et tout entier $p$. En divisant par $p$, en int\'egrant, et en passant \`a la limite, nous obtenons bien le r\'esultat voulu : la quantit\'e 
$$\tau_{B,\nu}(f)=\mathrm{lim}_{p\to \infty}\frac{1}{p}\int_{B}\Phi(\widetilde{df^{p}})\nu^{n}$$ ne d\'epend pas de la trivialisation choisie. L'application $\tau_{B,\nu} : \Gamma_{B,\nu} \to \R $ est le quasi-morphisme homog\`ene annonc\'e.

\subsection{Le quasi-morphisme $\mathfrak{S}$}
Notons d'abord que, puisque la classe $[\omega]$ n'est jamais nulle, notre hypoth\`ese force $c_{1}(V)\neq 0$. Cela exclut par exemple les vari\'et\'es symplectiques telles que les surfaces $K3$ ou les tores (on pourra consulter \cite{entov} pour la construction d'un quasi-morphisme homog\`ene sur le groupe $\widetilde{\ham}(V,\omega)$, pour les vari\'et\'es de premi\`ere classe de Chern nulle). Dans le cas o\`u $V$ est ou bien $\mathbb{CP}^{2}$ ou bien $\mathbb{CP}^{1}\times \mathbb{CP}^{1}$ muni du produit de la forme symplectique standard par elle-m\^eme, des r\'esultats de M.~Gromov \cite{gromov} assurent que le groupe fondamental de $\ham(V,\omega)$ est fini. Tout quasi-morphisme homog\`ene sur $\widetilde{\ham}(V,\omega)$ descend donc sur $\ham(V,\omega)$. Notre construction fournit donc un nouvel exemple de quasi-morphisme homog\`ene sur les groupes $\ham(\mathbb{CP}^{2},\omega_{0})$ et $\ham(\mathbb{CP}^{1}\times \mathbb{CP}^{1},\omega_{0}\times \omega_{0})$.

Fixons un $S^{1}$-fibr\'e principal $\pi : M\to V$ de classe d'Euler \'egale \`a $2c_{1}(V)$. Notant $X$ le champ de vecteurs sur $M$ engendr\'e par l'action du cercle, on peut trouver une $1$-forme $\alpha$ sur $M$ telle que $\alpha (X)=1$ et $d\alpha=\pi^{*}(s\omega)$. Nous sommes dans la situation du paragraphe $2.1$. Fixons \'egalement une structure presque-complexe $J$ sur $V$, compatible avec $\omega$. Le fibr\'e vectoriel $TV$ devient alors un fibr\'e hermitien, dont on peut choisir une trivialisation au-dessus d'un recouvrement $\{U_{\beta}\}$, avec des applications de transition $g_{\beta \gamma} : U_{\beta}\cap U_{\gamma}\to U(n)$, \`a valeurs dans le groupe des matrices unitaires de taille $n$. La famille d'applications $(det^{2}(g_{\beta \gamma}))$, d\'etermine un fibr\'e en cercles $E$ au-dessus de $V$, qui est isomorphe \`a $M$. Si l'on note $\Lambda(V)$ le fibr\'e en grassmannienne lagrangienne au-dessus de $V$, on a une application \mbox{$det^{2} : \Lambda(V)\to E$} qui n'est autre qu'une version fibr\'ee de l'application d\'ej\`a rencontr\'ee dans le cas lin\'eaire. Dans une trivialisation $U_{\gamma}\times \C^{n}$, un \'el\'ement $L$ de $\Lambda(V)$ s'\'ecrit $(x,u_{\gamma}(\RR))$, pour une matrice unitaire $u_{\gamma}$. On lui associe l'\'el\'ement $(x,det^{2}(u_{\gamma}))$ dans la trivialisation correspondante de $E$. En choisissant un isomorphisme entre $E$ et $M$ on obtient une application $\varphi : \Lambda(V)\to M$. Elle n'est bien s\^ur pas unique, le choix de $J$ et l'isomorphisme entre $E$ et $M$ interviennent. Cependant, elle a la vertu suivante : en restriction \`a chaque fibre, elle induit un isomorphisme entre les groupes fondamentaux. Une autre application $\varphi' : \Lambda(V)\to E$ construite par le m\^eme proc\'ed\'e serait donc de la forme $$\varphi'(L)=\chi(\pi(L))\cdot e^{2i\pi \kappa(L)}\cdot \varphi(L),$$ pour des application $\chi : V \to S^{1}$ et $\kappa : \Lambda(V)\to \R$.

Passons \`a la construction de notre dernier quasi-morphisme. On consid\`ere une isotopie hamiltonienne $(f_{t})$ engendr\'ee par le champ de vecteurs $Z_{t}$ (avec $\iota_{Z_{t}}\omega=dH_{t}$, $\int_{V}H_{t}\omega^{n}=0$). On note toujours $\Theta(f_{t})$ l'isotopie de $M$ engendr\'ee par le champ de vecteurs $\widehat{Z_{t}}-(H_{t}\circ \pi)X$. Si $L\in \Lambda(V)$, les deux courbes $$ \varphi(df_{t}\cdot L)\; \;  \mathrm{et}\; \;  \Theta(f_{t})(\varphi(L))$$ 
dans $M$, sont issues du m\^eme point et rel\`event la m\^eme courbe de $V$. Bien que le fibr\'e en cercles $M$ ne soit pas trivial, on peut se servir de la courbe $\Theta(f_{t})(\varphi(L))$ comme d'une horizontale ``le long du chemin $f_{t}(\pi(L))$" pour mesurer le nombre de rotation de la courbe $\varphi(df_{t}\cdot L)$. On peut donc \'ecrire $\varphi(df_{t}\cdot L)=e^{2i\pi \vartheta(t))}\cdot \Theta(f_{t})(\varphi(L))$ et d\'efinir une fonction continue sur $\Lambda(V)$ par $\ang(L,\{f_{t}\})=\vartheta(1)-\vartheta(0)$. Elle satisfait la relation : 
$$\ang (L,\{f_{t}\ast g_{t}f_{1}\}=\ang(L,\{f_{t}\})+\ang(df_{1}\cdot L, \{g_{t}\}).$$ 
\begin{proposition}
Pour toute paire de lagrangiens $(L_{0},L_{1})$ contenus dans la m\^eme fibre de $\Lambda(V)\to V$, et toute isotopie hamiltonienne $\{f_{t}\}$, nous avons : 
$$\vert \ang(L_{0},\{f_{t}\})-\ang(L_{1},\{f_{t}\})\vert \le 2n. $$
\end{proposition}
\textit{Preuve} : c'est une version fibr\'ee des r\'esultats du paragraphe $3.1$, qui permettent de construire le quasi-morphisme homog\`ene sur le rev\^etement universel du groupe symplectique. \hfill $\Box$
 
Nous d\'efinissons alors une fonction mesurable born\'ee sur $V$ par : $\ang(x,\{f_{t}\})=\mathrm{inf}_{L\in \Lambda(V)_{x}}\ang(L,\{f_{t}\})$. Elle satisfait : $$\vert \ang(x,\{f_{t}\ast g_{t}f_{1}\})-\ang(x,\{f_{t}\})-\ang(f_{1}(x),\{g_{t}\})\vert \le 6n.$$
L'application 
$$\begin{array}{rcl}
\widetilde{\ham}(V,\omega) & \to & \R \\
\{f_{t}\} & \mapsto & \int_{V}\ang(-,\{f_{t}\})\omega^{n} \\
\end{array}$$
est donc un quasi-morphisme. Si l'application $\varphi$ est modifi\'ee en une application $\varphi'$, comme expliqu\'e ci-dessus, la fonction $\ang$ se trouve chang\'ee en la fonction $\ang'(x,\{f_{t}\})$ \'egale \`a : 
$$\ang(L,\{f_{t}\})+\kappa (df_{1}\cdot L)-\kappa (L)+ \int_{0}^{1}\beta (X_{t})(f_{t}(x))dt,$$
si $\beta$ d\'esigne la $1$-forme ferm\'ee $d(\frac{ln \chi}{2i\pi})$. Nous avons donc : 
$$\vert \ang'(x,\{f_{t}\})-\ang(x,\{f_{t}\})-\int_{0}^{1}\beta(X_{t})(f_{t}(x))dt\vert \le 4n + 2 \mathrm{sup}_{\Lambda(V)}\vert \kappa\vert.$$ L'int\'egrale $\int_{V}\int_{0}^{1}\beta (X_{t})dt \omega^{n}$ \'etant nulle, on a : 
$$\vert \int_{V}\ang'(-,\{f_{t}\})\omega^{n}-\int_{V}\ang(-,\{f_{t}\})\omega^{n}\vert \le (4n + 2 \mathrm{sup}_{\Lambda(V)}\vert \kappa \vert)\mathrm{vol}(V).$$
L'homog\'en\'eis\'e 
$$\mathfrak{S}(\{f_{t}\})=\mathrm{lim}_{p\to \infty}\frac{1}{p}\int_{V}\ang(-,\{f_{t}\}^{p})\omega^{n}$$ de notre quasi-morphisme ne d\'epend donc pas du choix de l'application $\varphi$, et donc pas de $J$. Il ne d\'epend pas non plus du choix de la $1$-forme $\alpha$ (v\'erifiant $\alpha(X)=1$ et $d\alpha=\pi^{*}(s\omega)$). Nous calculons maintenant sa restriction sur les isotopies \`a support dans une boule $i : B \hookrightarrow V$. 

\begin{proposition}
Si l'isotopie $\{f_{t}\}$ est \`a support dans $B$, nous avons $$\mathfrak{S}(\{f_{t}\})=\tau_{B,\omega}(f_{1})+s\calabi_{B}(f_{1}).$$
\end{proposition}

\textit{Preuve} : on fixe une trivialisation unitaire du fibr\'e tangent au-dessus de $B$, et une trivialisation du fibr\'e en cercles $M$ au-dessus de $B$. On note $\lambda$ la primitive de $\omega$ sur $B$ telle que $\alpha =d(\frac{ln z}{2i\pi})+s\lambda$ dans cette trivialisation ($z$ d\'esigne la coordonn\'ee sur le cercle). L'application $\varphi$ lue dans cette trivialisation est de la forme :

$$\begin{array}{rcl}
B\times \Lambda(\R^{2n}) & \to & B \times S^{1} \\
(x,L)  & \mapsto & (x,e^{2i\pi \kappa(x,L)}\cdot det_{\RR}^{2}(L)) \\  
\end{array}$$ 

\noindent pour une application $\kappa : B\times \Lambda(\RR)\to \R$. On consid\`ere alors une isotopie hamiltonienne $(f_{t})$ \`a support contenu dans un compact $K$ de $B$, engendr\'ee par le champ de vecteurs $Z_{t}$ (avec $\iota_{Z_{t}}\omega=d\widetilde{H}_{t}$, o\`u $\widetilde{H}_{t}$ est d'int\'egrale nulle sur $V$, constante hors de $B$). Nous avons d'une part : 
$$\varphi(f_{t}(x),df_{t}\cdot L)=(f_{t}(x),e^{2i\pi \kappa(f_{t}(x),df_{t}\cdot L)}\cdot det^{2}_{\RR}(df_{t}\cdot L))$$
et d'autre part : 
$$\Theta(f_{t})(\varphi(x,L))=(f_{t}(x),e^{-2i\pi \int_{0}^{t}(\lambda(Z_{u})+s\widetilde{H}_{u})(f_{u}(x))du}\cdot e^{2i\pi \kappa(x,L)}\cdot det^{2}_{\RR}(L)).$$
La valeur de $\ang(L,\{f_{t}\})$ est donc : $$\Delta(det^{2}_{\R^{n}}(df_{t}(x)\cdot L))+\int_{0}^{1}(s\lambda(X_{t})+sH_{t})(f_{t}(x))dt+ \kappa(f_{1}(x),df_{1}(x)\cdot L)-\kappa(x,L).$$
Notons $C$ le maximum de $\vert \kappa \vert $ sur $\mathrm{supp}(f)\times \Lambda(\R^{2n})$. Nous obtenons 
$$\vert \ang(x,\{f_{t}\})-\Delta(det^{2}_{\R^{n}}(df_{t}(x)\cdot L))-s\int_{0}^{1}(\lambda(X_{t})+\widetilde{H}_{t})(f_{t}(x))dt\vert \le 2C+2n.$$ Hors de la boule $B$ la fonction $\ang(x,\{f_{t}\})$ est \'egale \`a $-s\int_{0}^{1}\widetilde{H}_{t}(f_{t}(x))dt$. En tenant compte de l'in\'egalit\'e $\vert \Delta(det^{2}_{\R^{n}}(df_{t}(x)\cdot L))-\Phi(\{df_{t}(x)\})\vert \le 2n$, nous obtenons : 

$$\vert \int_{V}\ang(-,\{f_{t}\}\nu^{n}-s\int_{0}^{1}\int_{B}\lambda(X_{t})dt\nu^{n}-\int_{B}\Phi(\{df_{t}(x)\})\nu^{n}\vert \le (2C+4n)\mathrm{vol}(B).$$

\noindent La m\^eme estimation reste vraie pour les it\'er\'es de $f$ car leur support est contenu dans celui de $f$. Nous obtenons donc : $\mathfrak{S}(\{f_{t}\})=s\int_{0}^{1}\int_{B}\lambda(X_{t})dt\omega^{n}+\tau_{B,\nu}(f_{1})=s\calabi (f_{1}) +\tau_{B,\nu}(f_{1})$. \hfill $\Box$

\medskip

\textbf{Remerciements.} \'Etienne Ghys m'a propos\'e de r\'efl\'echir \`a ce sujet et m'a encourag\'e; je le remercie pour tout cela.  Je tiens \'egalement \`a remercier Emmanuel Giroux et Bruno S\'evennec pour de nombreuses discussions. Enfin, je tiens \`a remercier vivement Leonid Polterovich qui m'a signal\'e une erreur dans une version ant\'erieure du th\'eor\`eme $2$.

\bibliographystyle{plain}
\bibliography{bibcal}

\medskip

\begin{flushleft}

\nopagebreak{ Pierre Py\\ Unit\'e de Math\'ematiques Pures et Appliqu\'ees \\
UMR 5669 CNRS \\
\'Ecole Normale Sup\'erieure de Lyon\\ 46, All\'ee d'Italie\\ 69364 Lyon 
Cedex 07 \\ FRANCE\\ Pierre.Py@umpa.ens-lyon.fr}
\end{flushleft}

\end{document}